\documentclass[opre,sglanonrev]{informs4}
\RequirePackage{tgtermes}
\RequirePackage{newtxtext}
\RequirePackage{newtxmath}
\RequirePackage{bm}
\RequirePackage{endnotes}

\OneAndAHalfSpacedXI % Current default line spacing
\usepackage{algorithm}
\usepackage{algpseudocode}
\usepackage{tikz}
\usepackage{natbib}
 \bibpunct[, ]{(}{)}{,}{a}{}{,}%
 \def\bibfont{\small}%
\EquationsNumberedThrough    % Default: (1), (2), ...
\TheoremsNumberedThrough     % Preferred (Theorem 1, Lemma 1, Theorem 2)
\ECRepeatTheorems  %  

\MANUSCRIPTNO{}

\usepackage{mathtools}
\usepackage[T1]{fontenc}
\usepackage{fix-cm}
\usepackage{enumitem}

\usepackage{caption}
\usepackage[labelfont=sf]{subcaption}
\def\E{{\mathbb E}}
\def\P{{\mathbb P}}

\def\cu{{\mathcal U}}
\def\cx{{\mathcal X}}

\usepackage[linktocpage,colorlinks,linkcolor=blue,anchorcolor=blue,citecolor=blue,urlcolor=blue,pagebackref]{hyperref}
\hypersetup{}
\renewcommand*{\backref}[1]{\ifx#1\relax \else ~ \fi}
\begin{document}
%%%%%%%%%%%%%%%%

%\setcounter{page}{1}
%\DOI{10.1287/xxxx.0000.0000}%

\RUNAUTHOR{Chen, Xin, and Zhao}

% \RUNTITLE{Predictive Maintenance in Manufacturing}
\RUNTITLE{Hidden Convexity in Queueing Models}

% \TITLE{Optimal Resource Allocation in Humanitarian Logistics: A Stochastic Programming Approach}
\TITLE{Hidden Convexity in Queueing Models}

\ARTICLEAUTHORS{%
\AUTHOR{Xin Chen}
\AFF{H. Milton Stewart School of Industrial and Systems Engineering, Georgia Tech, \EMAIL{xin.chen@isye.gatech.edu}}

\AUTHOR{Linwei Xin}
\AFF{School of Operations Research and Information Engineering, Cornell University, \EMAIL{lx267@cornell.edu}}

\AUTHOR{Minda Zhao}
\AFF{H. Milton Stewart School of Industrial and Systems Engineering, Georgia Tech, \EMAIL{mindazhao@gatech.edu}}
} % end of the block

\ABSTRACT{%
% Enter your abstract
We study the joint control of arrival and service rates in queueing systems with the objective of minimizing long-run expected cost minus revenue. Although the objective function is non-convex, first-order methods have been empirically observed to converge to globally optimal solutions. This paper provides a theoretical foundation for this empirical phenomenon by characterizing the optimization landscape and identifying a hidden convexity: the problem admits a convex reformulation after an appropriate change of variables. Leveraging this hidden convexity, we establish the Polyak-Łojasiewicz-Kurdyka (PŁK) condition for the original control problem, which excludes spurious local minima and supports global convergence guarantees for first-order methods. Our analysis applies to a broad class of $GI/GI/1$ queueing models, including those with Gamma-distributed interarrival and service times, as well as $GI/M/1$ queues with log-concave interarrival times. As a key ingredient in the proof, we establish a new convexity property of the expected queue length under a square-root transformation of the traffic intensity. 
}%

% \FUNDING{This research was supported by [grant number, funding agency].}

%Supplemental Material:
%Data Ethics & Reproducibility Note:

%\KEYWORDS{Stochastic programming, Decision support,Uncertainty, Disaster response, Optimization}

%\HISTORY{Received: Month DD, YYYY; Accepted: Month DD, YYYY; Published Online: Month DD, YYYY}
% \HISTORY{\today{}}

\maketitle
%%%%%%%%%%%%%%%%%%%%%%%%%%%%%%%%%%%%%%%%%%%%%%%%%%%%%%%%%%%%%%%%%%%%%%

% Text of your paper here

\section{Introduction}
This work studies a queueing model in which a service provider manages the system by jointly controlling the arrival and service rates, a topic that has been well studied in the literature over the past decades (e.g., \citealp{mendelson1985pricing, maglaras2003pricing, lee2014optimal, chen2024online}). We follow the standard setting where the service provider seeks the optimal arrival rate $\lambda^*$ and service rate $\mu^*$ that minimize the long-run expected cost minus revenue in a $GI/GI/1$ queue: $ \E [ Q_\infty (\lambda, \mu) ] + c(\mu) - r(\lambda)$, where $Q_\infty(\lambda, \mu)$ denotes the steady-state queue length under arrival rate $\lambda$ and service rate $\mu$; $c(\mu)$ is the cost of providing service rate $\mu$; and $r(\lambda)$ is the revenue corresponding to the arrival rate $\lambda$.

Solving this joint control problem over arrival and service rates is challenging because the objective function is non-convex \citep{harel1990convexity}. Existing studies often restrict attention to tractable settings such as $M/M/1$ queues for exact analysis \citep{ata2006dynamic} or asymptotic heavy-traffic regimes \citep{lee2014optimal, lee2019pricing}. More recently, \citet{chen2024online} propose a gradient-based method to numerically solve the problem. Surprisingly, despite the non-convex nature of the underlying problem, both \citet{chen2024online} and \citet{hu2024multi}, who control the arrival rate through pricing, observe that first-order methods empirically converge to global optimal solutions. This phenomenon reveals a striking gap between theory and practical performance. Understanding the mechanisms behind the success of first-order methods in these non-convex queueing models is essential for developing more principled online learning and offline optimization algorithms.

In this paper, we provide a theoretical foundation for this empirical phenomenon by identifying a hidden convexity. The key step is a reparameterization of the optimization problem through a change of variables. Under mild structural assumptions, we show that the original non-convex objective can be reformulated as a convex function over the new variables, subject to convex constraints. This convex reformulation applies to a broad class of $GI/GI/1$ queues, including $\Gamma_k/\Gamma_m/1$ queues whose interarrival and service times follow Gamma distributions with shape parameters $k$ and $m$ no less than $1$, as well as $GI/M/1$ queues with log-concave interarrival times. It also extends to Jackson networks of $M/M/1$ queues, $M/GI/1$ queues, and well-known approximate formulations of $GI/GI/s$ queues. Finally, leveraging this reformulation, we establish that the original problem satisfies the Polyak-Łojasiewicz-Kurdyka (PŁK) condition \citep{polyak1963gradient, lojasiewicz1963topological, kurdyka1998gradients}, which ensures the global convergence of first-order methods and helps explain the observed empirical success under standard assumptions.

A key ingredient in our analysis is a new convexity result. To construct the reformulation, we require the expected queue length to be convex after a square-root transformation of the traffic intensity $\rho$. Specifically, let $L_q(\rho)$ denote the expected queue length. We prove that $\hat{L}_q(\tau) \coloneqq L_q(\sqrt{\tau})$ is convex for $\Gamma_k/\Gamma_m/1$ queues with $k, m \ge 1$ and $GI/M/1$ queues with log-concave interarrival times. Interestingly, using Spitzer's identity, $\hat L_q(\tau)$ can be represented as the sum of a finite collection of potentially non-convex terms and infinitely many, rapidly vanishing convex terms. Because of this mixture, the global convexity of $\hat{L}_q(\tau)$ is far from trivial, which may be of independent interest.

\subsection{Literature Review}\label{subsection: literature}
In recent years, learning and optimization in queueing problems have attracted considerable attention (see, for example, \citealp{chen2023steady, chen2023incentivizing, kanoria2024blind, zhong2024learning}). One powerful tool for solving these problems is convexity, which has been studied for decades in queueing systems. Much of the literature establishes univariate convexity of performance measures while holding other parameters fixed, e.g., convexity in the service rate $\mu$ with fixed arrival rate $\lambda$ \citep{weber1983note}, and convexity in $\lambda$ with fixed $\mu$ \citep{fridgeirsdottir2005note}. Related results also show convexity with respect to the number of servers \citep{weber1980note}. By contrast, joint convexity in arrival and service rates generally fails. \citet{harel1990convexity} demonstrates that ``neither the expected number of customers in the queue nor in the system is jointly convex for the queues $M/M/x$, $M/G/1$."

Despite the loss of joint convexity, the optimization problem remains tractable when the decision-maker’s objective is solely to minimize queue performance metrics, such as the expected queue length. Existing work establishes convexity of various performance measures in traffic intensity $\rho$ for $M/M/c$ queues (see details in \citealp{harel1987strong}). The difficulty arises when service cost and revenue are involved, as the convexity fails due to the non-convex constraint $\rho = \lambda/\mu$. In addition, \citet{agrawal2020differentiating} apply a log-log convex programming approach to minimize the weighted sum of the service loads in $M/M/c$ queues. However, their formulation does not contain the service cost minus revenue term, which is generally not log-log convex. \citet{chen2024online} impose assumptions ensuring the joint strong convexity of the long-run expected cost. The assumptions are limited to $M/GI/1$ queues, which allow for a closed-form expression of the expected system length as a function of traffic intensity. Departing significantly from previous approaches, our work is the first to reveal the hidden convexity in the joint control of arrival and service rates for $GI/GI/1$ queues. Under mild conditions, we construct convex reformulations and establish the PŁK condition of the original problem. We contribute to the optimization literature by identifying a class of queueing control problems that satisfy the PŁK condition, thereby complementing recent advances on tractable non-convex optimization problems in operations and control models \citep{bhandari2024global, chen2024landscape}.

Our use of square-root transformations is conceptually related to, but technically distinct from, the classical line of work in queueing and stochastic networks. In capacity-allocation problems for Jackson networks, square-root allocation formulas have a long history, dating back to Kleinrock's work \citep{kleinrock1964communication} and subsequent developments for product-form and generalized Jackson networks \citep{wein1989capacity,dieker2017optimal}. These formulas typically arise from applying Lagrangian methods to approximately solve capacity-allocation problems. In other models, such as multiclass queueing networks, entropy-type formulations arise in the large-deviation objective \citep{pittel1979closed} and in the proportional-fairness interpretation of such networks \citep{walton2009proportional}. Our contribution is different from these classical analyses: we use square-root transformations as reparameterizations of the original decision variables to reveal hidden convexity and establish a P{\L}K landscape for the underlying nonconvex queueing-control problem.

Our paper is also related to the rich literature on joint pricing and capacity sizing in queueing systems, where both the arrival and service rates are controlled via a transformation between price and arrival rate. This transformation is commonly used in the revenue management literature \citep{talluri2005theory} to address the non-concavity of revenue with respect to price. Several existing studies focus on static pricing control, in which the price is independent of the queue length. Since a closed-form expression for the expected queue length is generally unavailable in $GI/GI/1$ queues, some studies adopt heavy-traffic approximations to derive asymptotically optimal arrival and service rates \citep{lee2014optimal, lee2019pricing}. Differently, \citet{chen2023online} and \citet{chen2024online} apply online first-order methods to solve the problem. In the dynamic setting, \citet{ata2006dynamic} investigate the optimal dynamic pricing and service-rate control policy in $M/M/1$ queues, while \citet{kim2018value} use an asymptotic large-system approach.

\section{Problem Setting}\label{sec-setting}
Consider a single-server $GI/GI/1$ queueing system with customer arrivals according to a stochastic process. The interarrival times $\{T_n\}_{n=1}^\infty$ are independent and identically distributed (i.i.d.) with a general distribution, corresponding to the first ``$GI$." Similarly, the service times $\{S_n\}_{n=1}^\infty$ are also i.i.d. with a general distribution, representing the second ``$GI$." The sequences $\{T_n\}_{n=1}^\infty$ and $\{S_n\}_{n=1}^\infty$ are mutually independent. We model the queueing system by allowing control over both the arrival rate $\lambda > 0$ and service rate $\mu > 0$, such that the interarrival and service times are given by $T_n / \lambda$ and $S_n / \mu$, respectively; namely, each interarrival (service) time is expressed as the product of the inverse deterministic arrival (service) rate and a random variable. This decomposition is a standard approach in the queueing literature on rate control (e.g., \citealp{lee2014optimal, whitt2015stabilizing, whitt2018time, chen2024online}). A natural interpretation is that $S_n$ may represent the size of a message to be transmitted in a communication network, while $\mu$ corresponds to the processing rate of messages. Without loss of generality, we assume $\E[T_n] = \E[S_n] = 1$, which can be achieved by appropriately rescaling the arrival and service rates.

By controlling the arrival and service rates, the system receives revenue $r(\lambda)$ and incurs a non-decreasing capacity cost $c(\mu)$ per time unit. The decision-maker searches for the optimal arrival rate $\lambda^*$ and service rate $\mu^*$ to minimize the long-run expected cost over a bounded domain:
\begin{equation}\label{opt-main problem}
    \begin{aligned}
        \min_{\lambda, \mu} \quad & \E \left[ Q_\infty (\lambda, \mu) \right] + c(\mu) - r(\lambda)\\
        \text{s.t.} \quad & \underline{\lambda} \le \lambda \le \bar{\lambda}, \quad \underline{\mu} \le \mu \le \bar{\mu},
    \end{aligned}
\end{equation}
with $\underline{\lambda} > 0$. To ensure stability of the queueing system and to state the landscape results in a differentiable form, we impose the following standing assumptions.

\begin{assumption}\label{assumption-general model}
    All of the following conditions hold:
    \begin{enumerate}
        \item\label{assumption-differentiable} The revenue $r(\lambda)$ and service cost $c(\mu)$ are continuously differentiable.
        \item The bounds $\bar{\lambda}$ and $\underline{\mu}$ satisfy $\bar{\lambda} < \underline{\mu}$, implying that the service system is uniformly stable.
        \item The steady-state expected queue length is finite and differentiable on the feasible region.
    \end{enumerate}
\end{assumption}

The first two conditions are adopted from \citet{chen2024online}. Under uniform stability, finite second moments are standard sufficient conditions for the finiteness of the steady-state expected queue length. We impose the differentiability condition because the P{\L}K condition is formulated in terms of the gradient of the objective. This differentiability condition can be verified directly for the distributional classes analyzed later, including the exponential and Gamma cases.

In general, \eqref{opt-main problem} is non-convex \citep{harel1990convexity}. To address this issue, we introduce two different approaches to construct convex reformulations of \eqref{opt-main problem} in the following section. The analysis applies to a broad class of $GI/GI/1$ queueing systems, including $\Gamma_k/\Gamma_m/1$ queues with $k,m \ge 1$ and $GI/M/1$ queues with log-concave interarrival times.

\section{Convex Reformulation}\label{sec-reformulation}
We begin by analyzing the expected queue length $\E [ Q_\infty (\lambda, \mu) ]$, whose closed-form expression is generally unavailable for $GI/GI/1$ queues. Nevertheless, we can relate it to the expected waiting time via Little’s law \citep{little1961proof} by $\E [ Q_\infty (\lambda, \mu) ] = \lambda \E [W_\infty (\lambda, \mu)]$. Here, $W_\infty (\lambda, \mu)$ denotes the stationary waiting time, i.e., the limit of waiting times $W_n(\lambda, \mu)$ as $n\to\infty$ and the sequence $\{W_n(\lambda, \mu)\}_n$ satisfies Lindley's equation \citep{lindley1952theory}:
\begin{equation*}
    W_{n+1}(\lambda, \mu) = \left( W_n(\lambda, \mu) + \frac{S_n}{\mu} - \frac{T_n}{\lambda} \right)^+.
\end{equation*}
Under Assumption~\ref{assumption-general model}, the increments have negative drift with $\E[S_i / \mu - T_i / \lambda] < 0$. This ensures the stability and allows us to apply the celebrated Spitzer's identity (e.g., \citealp{kingman1962some}):
\begin{equation*}
    \E \left[ W_{\infty}(\lambda, \mu) \right] = \sum_{n=1}^\infty \frac{1}{n} \E \left[ \left(\sum_{i=1}^n \frac{S_i}{\mu} - \sum_{i=1}^n \frac{T_i}{\lambda} \right)^+ \right].
\end{equation*}
Define
\begin{equation*}
    l_n(\rho) \coloneqq \frac{1}{n} \E \left[ \left(\sum_{i=1}^n \rho S_i - \sum_{i=1}^n T_i \right)^+ \right] \quad \text{and} \quad L_q(\rho) \coloneqq \sum_{n=1}^\infty l_n(\rho).
\end{equation*}
For general cases, $L_q(\rho)$ is defined on the stability interval $\rho < \mathbb E[T_1]/\mathbb E[S_1]$. In the unit-mean normalization, the stability interval reduces to $\rho\in(0,1)$. By Little's law, it follows that $\E [ Q_\infty (\lambda, \mu) ] = L_q(\lambda / \mu)$. For any realization of $\{S_i, T_i\}$, the function $(\sum_{i=1}^n \rho S_i - \sum_{i=1}^n T_i)^+$ is a maximum of two linear functions and is therefore convex. Since convexity is preserved under expectation and countable summation, $L_q(\rho)$ is convex. However, the composition term $\lambda/\mu$ destroys the joint convexity. We present two approaches for constructing a convex reformulation: the first involves a change of variables with respect to $\lambda$, and the second involves a change of variables with respect to $\mu$.

\subsection{Reformulation I}\label{section: reformulation I}
Define $\hat{\lambda} \coloneqq \sqrt{\lambda}$, let $\pi(x,y)\coloneqq x^2/y$ for $y>0$, and define $\hat{r}(\hat{\lambda}) \coloneqq r(\hat{\lambda}^2)$. Then, (\ref{opt-main problem}) can be reformulated as
\begin{equation}\label{opt-gg1-r1}
    \tag*{(R1)}
    \begin{aligned}
        \min_{\hat\lambda, \mu} \quad & L_q \left( \pi (\hat\lambda, \mu) \right) + c(\mu) - \hat r(\hat\lambda)\\
        \text{s.t.} \quad & \sqrt{\underline\lambda} \le \hat\lambda \le \sqrt{\bar\lambda}, \quad \underline{\mu} \le \mu \le \bar{\mu}.
    \end{aligned}
\end{equation}
To ensure the convexity (or strong convexity) of \ref{opt-gg1-r1}, we need to make the following assumption.

\begin{assumption}\label{assumption:lambda_hat} 
    We impose one of the following two conditions:
    \begin{henumerate}
        \item \label{assumption:lambda_hat convex} $\hat{r}(\hat{\lambda})$ is concave in $\hat{\lambda}$ and $c(\mu)$ is convex in $\mu$.
        \item \label{assumption:lambda_hat strongly convex} $\hat{r}(\hat{\lambda})$ is $\kappa_1$-strongly concave in $\hat{\lambda}$ and $c(\mu)$ is $\kappa_1$-strongly convex in $\mu$.
    \end{henumerate}
\end{assumption}

\begin{theorem}\label{theorem-gg1-r1}
    Suppose Assumption~\ref{assumption-general model} holds. Under Assumption~\ref{assumption:lambda_hat}.\ref{assumption:lambda_hat convex}, the optimization problem~\ref{opt-gg1-r1} is convex. If, moreover, Assumption~\ref{assumption:lambda_hat}.\ref{assumption:lambda_hat strongly convex} holds, then problem~\ref{opt-gg1-r1} is $\kappa_1$-strongly convex.
\end{theorem}

\begin{proof}{Proof of Theorem \ref{theorem-gg1-r1}.}
    For any fixed realization of $\{S_i,T_i\}$, the map $\rho \mapsto \left(\sum_{i=1}^n \rho S_i - \sum_{i=1}^n T_i\right)^+$ is convex and nondecreasing in $\rho$. Therefore, $L_q(\rho)$ is convex and nondecreasing on $(0,1)$, since both properties are preserved under expectation and countable summation. Next, $\pi(\hat{\lambda},\mu)$ is jointly convex in $(\hat{\lambda},\mu)$, because it is the perspective of the convex quadratic function: $\pi(\hat{\lambda},\mu)=\mu (\hat{\lambda} / \mu )^2$. By Assumption~\ref{assumption-general model}, $\pi(\hat{\lambda},\mu)\in(0,1)$. Since $L_q$ is convex and nondecreasing on $(0,1)$, the composition rule for convex functions implies that $L_q(\pi(\hat{\lambda},\mu))$ is jointly convex in $(\hat{\lambda},\mu)$. 
    
    By Assumption~\ref{assumption:lambda_hat}.\ref{assumption:lambda_hat convex}, $c(\mu)-\hat r(\hat{\lambda})$ is jointly convex in $(\hat{\lambda},\mu)$. Hence, problem \ref{opt-gg1-r1} is convex. Moreover, by Assumption~\ref{assumption:lambda_hat}.\ref{assumption:lambda_hat strongly convex}, $c(\mu)-\hat r(\hat{\lambda})$ is jointly $\kappa_1$-strongly convex in $(\hat{\lambda},\mu)$. Adding the convex term $L_q(\pi(\hat{\lambda},\mu))$ preserves strong convexity. Therefore, problem \ref{opt-gg1-r1} is $\kappa_1$-strongly convex. \Halmos
\end{proof}

\begin{remark}
    Besides convexity, \ref{opt-gg1-r1} is second-order cone (SOC) representable for $M/M/1$ queue if $c(\mu)$ and $-\hat{r}(\hat{\lambda})$ are SOC representable. When $S_n, T_n\sim\exp(1)$, $L_q(\rho) = \rho^2/(1-\rho)$ and \ref{opt-gg1-r1} is equivalent to
    \begin{equation*}
        \begin{aligned}
            \min_{\hat\lambda, \mu} \quad & t_q + t_c + t_r\\
            \text{s.t.} \quad & \rho^2 / (1-\rho) \le t_q, \quad c(\mu) \le t_c, \quad -\hat r(\hat\lambda) \le t_r, \quad \rho \ge \hat\lambda^2 / \mu, \quad 0 \le \rho < 1, \\
            & \sqrt{\underline\lambda} \le \hat\lambda \le \sqrt{\bar\lambda}, \quad \underline{\mu} \le \mu \le \bar{\mu}.
        \end{aligned}
    \end{equation*}
    The constraint $\rho^2 / (1-\rho) \le t_q$ is equivalent to $2t_q(1-\rho) \ge \|\sqrt{2}\rho\|_2^2$ and $\rho \ge \hat\lambda^2 / \mu$ is equivalent to $2\rho\mu \ge \|\sqrt{2}{\hat\lambda}\|_2^2$. All constraints are SOC representable. Thus, the optimization problem is a second-order cone program (SOCP), which can be efficiently solved by solvers such as CPLEX and Gurobi.
\end{remark}

To interpret Assumption~\ref{assumption:lambda_hat}, we adopt the standard revenue model in which the price $p$ and the arrival rate $\lambda$ are linked by an invertible mapping. Under this specification, revenue can be written equivalently as $r(\lambda)=\lambda p(\lambda)$ or $r(p)=p\lambda(p)$. When $r$ is twice differentiable, $\hat r(\hat\lambda)=r(\hat\lambda^2)$ is concave if and only if $2r'(\lambda)+4\lambda r''(\lambda)\le 0$, and it is $\kappa_1$-strongly concave if $2r'(\lambda)+4\lambda r''(\lambda)\le -\kappa_1$ throughout the feasible region. As an example, consider the demand function $\lambda(p)=p^{-\gamma}$ with $\gamma\in(1,2]$. In this case, $\hat r(\hat\lambda)=\hat\lambda^{2-2/\gamma}$, which is concave for all $\gamma\in(1,2]$, and is strongly concave on any compact feasible interval when $\gamma\in(1,2)$. 

However, for many demand functions studied in revenue management \citep{talluri2005theory, besbes2015surprising}, Assumption~\ref{assumption:lambda_hat} may not hold globally.
\begin{hitemize}
    \item For linear demand $\lambda(p)=a-bp$ with $a,b>0$, we have $r(\lambda)=(a\lambda-\lambda^2)/b$, so that $2r'(\lambda)+4\lambda r''(\lambda)=(2a-12\lambda)/b$ and hence the condition holds whenever $\lambda\ge a/6$.
    
    \item For exponential demand $\lambda(p)=ae^{-bp}$ with $a,b>0$, we have $r(\lambda)=\frac{\lambda}{b}\log\!\left(\frac{a}{\lambda}\right)$. Thus, $2r'(\lambda)+4\lambda r''(\lambda)=2(\log(a/\lambda) - 3) / b$, which is nonpositive whenever $\lambda\ge ae^{-3}$.
    
    \item For logit demand $\lambda(p)\coloneqq \frac{\exp(a-p)}{1+\exp(a-p)}$, we have $\hat r(\hat \lambda)=\hat \lambda^2 \ln(1-\hat \lambda^2)-\hat \lambda^2\ln(\hat \lambda^2)+a \hat \lambda^2$, which is generally non-concave, although one can still characterize local regions on which concavity holds. For example, when $a=1$, $\hat r(\hat\lambda)$ is concave for $\hat\lambda>\hat\lambda^*/2$, where $\hat\lambda^*$ denotes the global maximizer of $\hat r(\hat\lambda)$.
\end{hitemize}
Hence, to fully accommodate these demand functions, we propose an alternative approach to constructing the convex reformulation.

\subsection{Reformulation II}\label{section: reformulation II}
Define $\hat{\mu}\coloneqq \mu^2$ and $\hat c(\hat\mu) \coloneqq c(\sqrt{\hat\mu})$. Then, (\ref{opt-main problem}) can be reformulated as
\begin{equation}\label{opt-gg1-r2}
\tag*{(R2)}
    \begin{aligned}
        \min_{\lambda, \hat{\mu}} \quad & L_q \left( \sqrt{\pi \left( \lambda, \hat\mu \right)} \right) + \hat c(\hat\mu) - r(\lambda)\\
        \text{s.t.} \quad & \underline{\lambda} \le \lambda \le \bar{\lambda}, \quad \underline{\mu}^2 \le \hat{\mu} \le \bar{\mu}^2.
    \end{aligned}
\end{equation}
Similar to \ref{opt-gg1-r1}, we require the following assumption to ensure the convexity (or strong convexity) of \ref{opt-gg1-r2}.

\begin{assumption} \label{assumption:mu_hat} 
    We impose one of the following two conditions:
    \begin{henumerate}
        \item \label{assumption:mu_hat convex} $\hat c(\hat\mu)$ is convex in $\hat\mu$ and $r(\lambda)$ is concave in $\lambda$.
        \item \label{assumption:mu_hat strongly convex} $\hat c(\hat\mu)$ is $\kappa_2$-strongly convex in $\hat\mu$ and $r(\lambda)$ is $\kappa_2$-strongly concave in $\lambda$.
    \end{henumerate}
\end{assumption}

When $c$ is twice differentiable, $\hat{c}(\cdot)$ is convex if and only if $\mu c''(\mu) - c'(\mu) \ge 0$. For example, $c(\mu) = \alpha \mu^\gamma$ with $\gamma \ge 2$ and $\alpha \ge 0$ is a simple class that satisfies Assumption~\ref{assumption:mu_hat}. Setting $\gamma = 2$ recovers the quadratic cost function commonly used in the literature \citep{ata2006dynamic, chen2024online}. Moreover, the strong concavity requirement of $r$ in Assumption~\ref{assumption:mu_hat}.\ref{assumption:mu_hat strongly convex} holds true for a wide range of demand functions. For example, under linear demand $\lambda(p)=a-bp$ with $a,b>0$, the revenue function $r(\lambda)=(a\lambda-\lambda^2)/b$ is a $2/b$-strongly concave quadratic. Under the logit demand model, $r(\lambda)=\lambda\ln(1-\lambda)-\lambda\ln(\lambda)+a\lambda$, which is $27/4$-strongly concave on $(0,1)$. Hence, \ref{opt-gg1-r2} is broadly applicable as long as the service cost remains convex under the square root transformation, which typically holds under mild conditions.

\begin{figure}[!htbp]
    \FIGURE
    {% Subcaptions and Graphics
        \subcaptionbox{Plot of $\hat{l}_n(\tau)$ for different $n$.\label{fig:l_n}}
        {\includegraphics[width=0.485\textwidth]{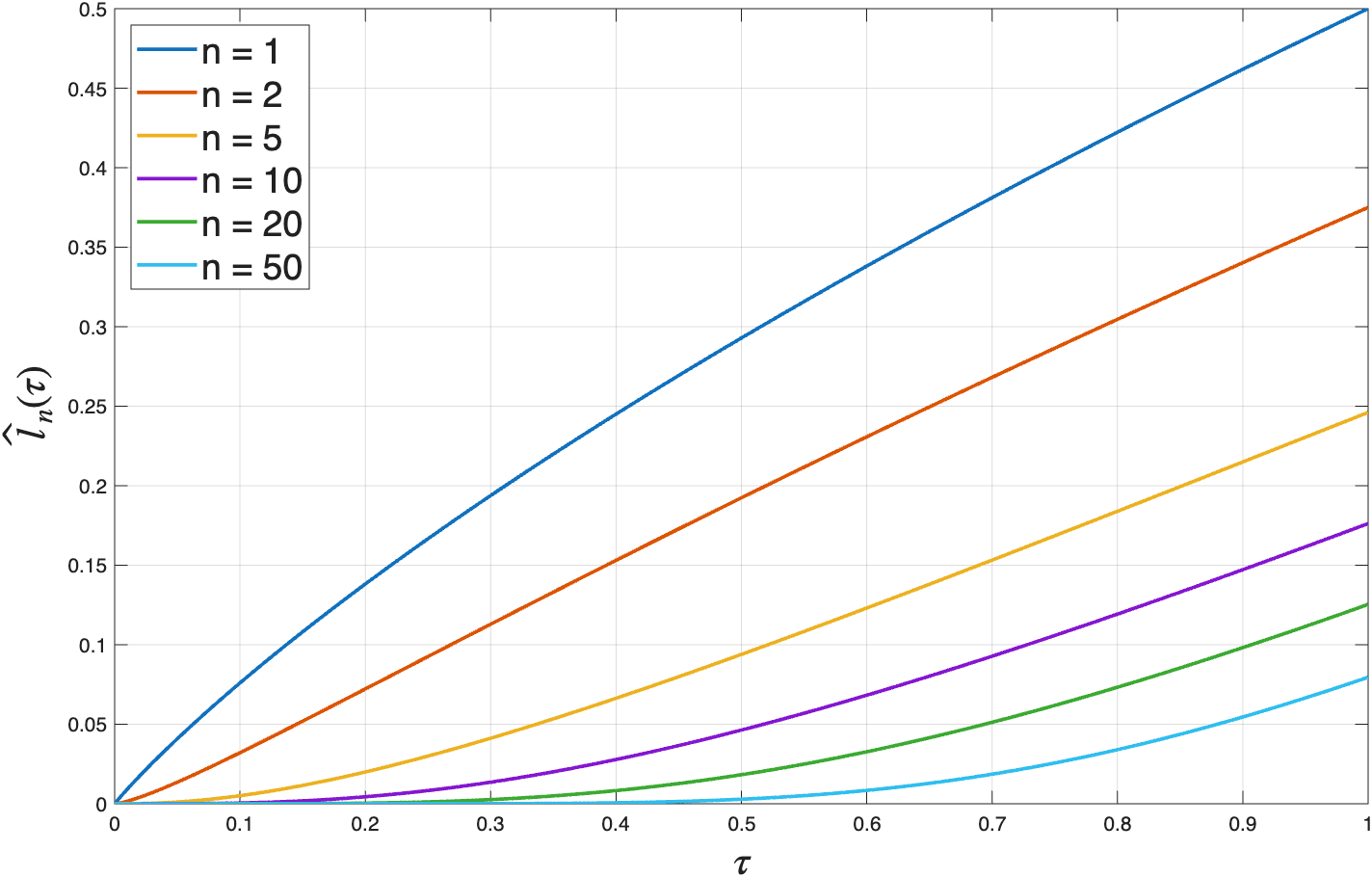}}
        \hfill\subcaptionbox{Plot of $\sum_{n=1}^N \hat{l}_n(\tau)$ for different $N$.\label{fig:sum l_n}}
        {\includegraphics[width=0.48\textwidth]{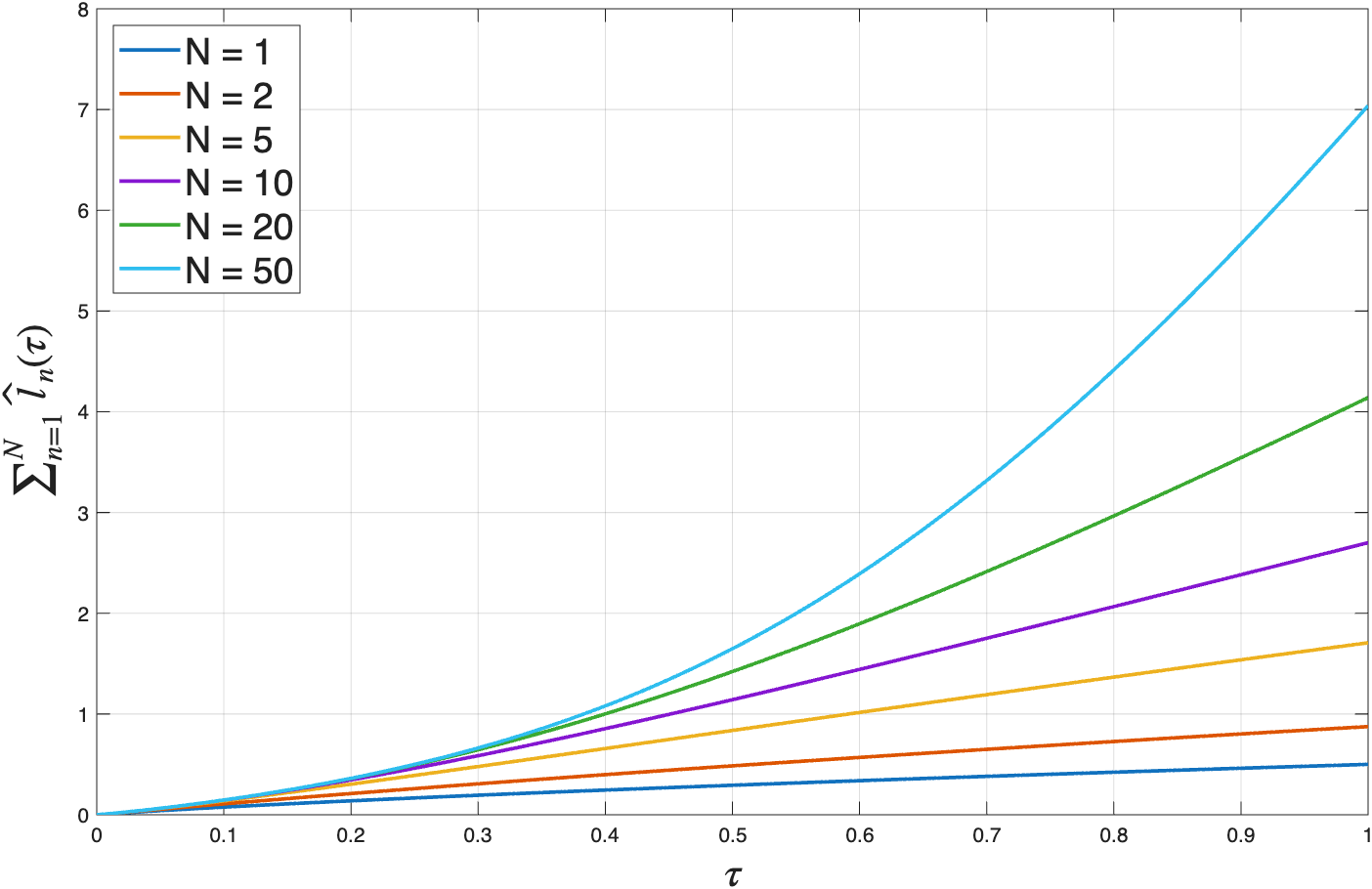}}
    }
    {% Main caption
    Plots of $\hat{l}_n(\tau)$ and $\sum_{n=1}^N \hat{l}_n(\tau)$ when $S_n, T_n \sim \exp(1)$.
    \label{fig:mm1}}
    {}
\end{figure}

By Assumption~\ref{assumption:mu_hat}, $\hat c(\hat\mu) - r(\lambda)$ is jointly convex in $(\lambda, \hat{\mu})$. For $\tau\in(0,1)$, define
\begin{equation*}
    \hat{l}_n(\tau) \coloneqq \frac{1}{n} \E \left[ \left(\sum_{i=1}^n \sqrt{\tau} S_i - \sum_{i=1}^n T_i \right)^+ \right] \quad \text{and} \quad \hat{L}_q(\tau) \coloneqq \sum_{n=1}^\infty \hat{l}_n(\tau).
\end{equation*}
Then, $\hat{L}_q ( \pi(\lambda, \hat{\mu}) ) = L_q ( \sqrt{\pi(\lambda, \hat\mu)} )$. To establish the convexity of \ref{opt-gg1-r2}, it is sufficient to prove that $\hat{L}_q(\tau)$ is convex because $\hat{L}_q(\tau)$ is non-decreasing and $\pi(\lambda, \hat\mu)$ is jointly convex. However, this is far from trivial, because $\hat{l}_n(\tau)$ can be non-convex (or even concave) for small $n$. Although $\hat{l}_n(\tau)$ eventually becomes convex as $n$ grows large, its magnitude seems to decay rapidly (see Figure~\ref{fig:l_n} for an illustration). Hence, it remains unclear whether the curvature of $\hat{L}_q(\tau)$ is primarily driven by the finitely many non-convex initial terms or by the infinitely many rapidly vanishing convex tails. For the $M/M/1$ queue, the summation $\sum_{n=1}^N \hat{l}_n(\tau)$ appears to be dominated by the convex terms $\hat{l}_n(\tau)$ corresponding to large $n$ as $N\to\infty$ (see Figure~\ref{fig:sum l_n}). In what follows, we rigorously establish the convexity result not only for the $M/M/1$ queue but also for $\Gamma_k/\Gamma_m/1$ queues with $k, m \ge 1$ and $GI/M/1$ queues with log-concave interarrival times.

\subsubsection*{Special Case: $M/M/1$ Queue}\label{subsec-mm1}
To build intuition, we start with the $M/M/1$ queue where $S_n, T_n\sim\exp(1)$. Then, the interarrival time $T_n / \lambda$ and the service time $S_n/\mu$ follow $\exp(\lambda)$ and $\exp(\mu)$, respectively. As a special case of the $GI/GI/1$ queue, the $M/M/1$ queue admits a closed-form expression of the expected queue length, e.g., $\E[Q_\infty(\lambda, \mu)] = \lambda^2 / [\mu(\mu - \lambda)]$. Hence, $\hat{L}_q(\tau)$ simplifies to $\tau / (1-\sqrt\tau)$.

\begin{theorem}\label{theorem-mm1}
    Suppose Assumption~\ref{assumption-general model} holds. In the case $T_n,S_n \overset{\text{i.i.d.}}{\sim} \exp(1)$, the function $\hat L_q(\tau)$ is convex on $(0,1)$. Under Assumption~\ref{assumption:mu_hat}.\ref{assumption:mu_hat convex}, the reformulated problem \ref{opt-gg1-r2} is convex. If, moreover, Assumption~\ref{assumption:mu_hat}.\ref{assumption:mu_hat strongly convex} is satisfied, then problem \ref{opt-gg1-r2} is $\kappa_2$-strongly convex.
\end{theorem}

\begin{proof}{Proof of Theorem \ref{theorem-mm1}.}
    For any $\tau \in (0,1)$, the geometric series expansion gives $\hat L_q(\tau)=\tau/(1-\sqrt{\tau})=\sum_{n=0}^\infty \tau^{1+n/2}$. Since each term $\tau^{1+n/2}$ is convex and increasing on $(0,1)$, it follows that $\hat L_q$ is convex and increasing on $(0,1)$. By Assumption~\ref{assumption-general model}, $\pi(\lambda,\hat\mu)\in(0,1)$. Since $\pi(\lambda,\hat\mu)$ is convex and $\hat L_q$ is convex and nondecreasing, the composition rule for convex functions implies that $\hat L_q(\pi(\lambda,\hat\mu))$ is jointly convex in $(\lambda,\hat\mu)$.

    Under Assumption~\ref{assumption:mu_hat}.\ref{assumption:mu_hat convex}, the function $\hat c(\hat\mu)-r(\lambda)$ is jointly convex in $(\lambda,\hat\mu)$. Hence, problem \ref{opt-gg1-r2} is convex. Moreover, by Assumption~\ref{assumption:mu_hat}.\ref{assumption:mu_hat strongly convex}, $\hat c(\hat\mu)-r(\lambda)$ is jointly $\kappa_2$-strongly convex in $(\lambda,\hat\mu)$. Since adding a convex function preserves strong convexity, problem \ref{opt-gg1-r2} is $\kappa_2$-strongly convex.\Halmos
\end{proof}

Although Theorem~\ref{theorem-mm1} is subsumed by more general Gamma distribution results presented in the following subsection, we state the theorem explicitly because the convex reformulation for $M/M/1$ queues extends naturally to other settings, including Jackson networks of $M/M/1$ queues, $M/GI/1$ queues, and approximate formulations of $GI/GI/s$ queues.

\subsubsection*{Extension: Jackson Network}
Consider a Jackson queueing network with $N$ stations. Each station $i \in [N]$ has an external Poisson arrival process with rate $\Lambda_i$. Revenue is collected upon admission of exogenous arrivals to the network. The service time at station $i$ follows an exponential distribution with rate $\mu_i$. A customer completing service at station $i$ will either move to some new station $j$ with probability $P_{ij}$ or leave the system with probability $1 - \sum_{j = 1}^N P_{ij}$. We assume the network is open, i.e., every customer can only visit finitely many stations before leaving the network with probability one. For example, the network is open if the strict inequality $\sum_{j = 1}^N P_{ij} < 1$ holds for any station $i\in[N]$. Then, the actual arrival rate $\lambda_i$ at station $i$ satisfies the balance equation
\begin{equation*}
    \lambda_i = \Lambda_i + \sum_{j=1}^N P_{ji} \lambda_j , \quad \forall i\in[N].
\end{equation*}

Applying Jackson's Theorem \citep{jackson1957networks, jackson1963jobshop}, the expected queue length at station $i$ is $\lambda^2_i/ [\mu_i(\mu_i - \lambda_i)]$. Then, we aim to solve the following optimization problem:
\begin{equation} \label{opt-Jackson}
    \begin{aligned}
        \min_{\lambda, \mu, \Lambda} \quad & \sum_{i = 1}^N \left[ \frac{\lambda^2_i}{\mu_i(\mu_i - \lambda_i)} + c(\mu_i) - r(\Lambda_i) \right]\\
        \text{s.t.} \quad & \lambda_i = \Lambda_i + \sum_{j=1}^N P_{ji} \lambda_j, &\forall i\in[N],\\
        & \underline{\lambda} \le \lambda_i \le \bar{\lambda}, \quad \underline{\mu} \le \mu_i \le \bar{\mu}, \quad 0 \le \Lambda_i \le \bar{\Lambda}, \quad  &\forall i\in[N].
    \end{aligned}
\end{equation}

Define $\hat{\mu}_i \coloneqq \mu_i^2$. Then, (\ref{opt-Jackson}) can be reformulated as
\begin{equation}\label{opt-Jackson reformulation}
    \begin{aligned}
        \min_{\lambda, \hat{\mu}, \Lambda} \quad & \sum_{i = 1}^N \left[ \frac{\lambda^2_i}{\sqrt{\hat\mu_i} \left( \sqrt{\hat\mu_i} - \lambda_i \right)} + \hat c(\hat{\mu}_i) - r(\Lambda_i) \right]\\
        \text{s.t.} \quad & \lambda_i = \Lambda_i + \sum_{j=1}^N P_{ji} \lambda_j, &\forall i\in[N],\\
        & \underline{\lambda} \le \lambda_i \le \bar{\lambda}, \quad \underline{\mu}^2 \le \hat{\mu}_i \le \bar{\mu}^2, \quad 0 \le \Lambda_i \le \bar{\Lambda}, \quad  &\forall i\in[N].
    \end{aligned}
\end{equation}

Under Assumption~\ref{assumption:mu_hat}.\ref{assumption:mu_hat convex}, each transformed expected queueing length is convex in $(\lambda_i,\hat\mu_i)$, each transformed service-cost term is convex in $\hat\mu_i$, and each revenue term $r(\Lambda_i)$ is concave in $\Lambda_i$. Since the balance equations and box constraints are affine, the reformulated Jackson-network problem (\ref{opt-Jackson reformulation}) is convex. Crucially, this exact network extension relies on the product-form stationary distribution of Jackson networks, which makes the steady-state queueing cost separable across stations. Without such separability, e.g., a generalized Jackson network, the single-station convexity argument does not directly lift to the network level.

\subsubsection*{Extension: $M/GI/1$ Queue}
Consider the setting where the interarrival times are exponential random variables, and the service times follow general distributions with finite variance $\text{Var}(S)<+\infty$. Applying the Pollaczek-Khinchine formula, \ref{opt-gg1-r2} reduces to
\begin{equation} \label{opt-mg1 formulation}
    \begin{aligned}
        \min_{\lambda, \hat\mu} \quad & \frac{1 + \text{Var}(S)}{2} \cdot \frac{\lambda^2}{\sqrt{\hat{\mu}} \left( \sqrt{\hat{\mu}} - \lambda \right)} + \hat c(\hat\mu) - r(\lambda)\\
        \text{s.t.} \quad & \underline{\lambda} \le \lambda \le \bar{\lambda}, \quad \underline{\mu}^2 \le \hat\mu \le \bar{\mu}^2.
    \end{aligned}
\end{equation}
The only difference between the optimization of $M/GI/1$ and $M/M/1$ queues lies in the coefficient in the objective. Therefore, we can similarly prove the convexity of (\ref{opt-mg1 formulation}).

\subsubsection*{Extension: $GI/GI/s$ Approximation}
Consider the setting where interarrival and service times follow general distributions with finite variance Var($T$) and Var($S$), respectively. Based on the approximation of $GI/GI/s$ system in \citet{shanthikumar2007queueing}, \ref{opt-gg1-r2} can be approximated by
\begin{equation}\label{opt-ggm reformulation}
    \begin{aligned}
        \min_{\lambda, \hat{\mu}} \quad & \frac{\text{Var}(T) + \text{Var}(S)}{2} \cdot \frac{ \left( \lambda / \sqrt{\hat{\mu}} \right)^{\sqrt{2(s+1)}}}{1-\lambda / \sqrt{\hat{\mu}}} + \hat c(\hat{\mu}) - r(\lambda)\\
        \text{s.t.} \quad & \underline{\lambda} \le \lambda \le \bar{\lambda}, \quad \underline{\mu}^2 \le \hat{\mu} \le \bar{\mu}^2.
    \end{aligned}
\end{equation}
Similar to Theorem~\ref{theorem-mm1}, the first term of the objective is a composition of two functions $\hat{L}_q (\pi(\lambda, \hat{\mu}) )$ with $\hat{L}_q (\tau) = \tau^{\sqrt{2(s+1)} / 2} / (1-\sqrt{\tau})$. For $s \ge 1$, $\hat{L}_q$ is monotone non-decreasing. Then, it is sufficient to check the convexity of $\hat{L}_q (\tau)$ by decomposing it via the Taylor's expansion for $\tau \in (0,1)$:
\begin{equation*}
    \begin{aligned}
        \hat{L}_q (\tau) = \frac{ \tau^{\frac{\sqrt{2(s+1)}}{2}}}{1-\sqrt{\tau}} = \tau^{\frac{\sqrt{2(s+1)}}{2}} \sum_{n=0}^\infty \left(\sqrt{\tau}\right)^n = \sum_{n=0}^\infty \tau^{\frac{\sqrt{2(s+1)} + n}{2}}.
    \end{aligned}
\end{equation*}
For each $n \ge 0$, the summand is convex, which implies that $\hat{L}_q (\tau)$ is convex. By the preservation under the composition of convex functions, we can establish the convexity of (\ref{opt-ggm reformulation}).

\subsubsection*{$\Gamma_k/\Gamma_m/1$ Queue}\label{subsec-gamma}
We turn to the more general setting of $\Gamma_k/\Gamma_m/1$ queues. To enforce $\E[T_n] = \E[S_n] = 1$, we set $T_n\sim\Gamma(k, k)$ and $S_n\sim\Gamma(m, m)$ under the (shape, rate) parameterization. Leveraging the structural properties of Gamma distributions, we establish the main result of this paper: $\hat{L}_q(\tau)$ is convex, despite being composed of finitely many non-convex components and infinitely many rapidly diminishing convex tails.

\begin{theorem}\label{theorem-ee1}
    Suppose that Assumption~\ref{assumption-general model} holds. In the case $T_n\sim\Gamma(k, k), S_n\sim\Gamma(m, m)$ with $k, m \ge 1$, the function $\hat{L}_q(\tau)$ is convex on $(0, 1)$. Under Assumption~\ref{assumption:mu_hat}.\ref{assumption:mu_hat convex}, the reformulated optimization problem \ref{opt-gg1-r2} is convex. If moreover, Assumption~\ref{assumption:mu_hat}.\ref{assumption:mu_hat strongly convex} is satisfied, then problem \ref{opt-gg1-r2} is $\kappa_2$-strongly convex.
\end{theorem}

Note that Theorem~\ref{theorem-ee1} holds for any $k, m \ge 1$. In the special case where both $k$ and $m$ are integers, we recover the Erlang queueing model $E_k/E_m/1$, a cornerstone in queueing theory.\vspace{-3pt}

\subsubsection*{$GI/M/1$ Queue with Log-concave Interarrival Times} 
In addition to the Gamma family with shape parameters no less than $1$, the convexity results also apply to a broad class of $GI/M/1$ queues. The next theorem shows that when service times are exponential and the interarrival-time density is log-concave on its support, the transformed expected queue length $\hat L_q(\tau)$ remains convex.

\begin{theorem}\label{theorem-GM1-logconcave}
    Suppose that Assumption~\ref{assumption-general model} holds. If $T_n$ has a log-concave density on its support, where the support is an interval contained in $[0,\infty)$, and $S_n \sim \exp(1)$, the function $\hat{L}_q(\tau)$ is convex on $(0, 1)$. Under Assumption~\ref{assumption:mu_hat}.\ref{assumption:mu_hat convex}, the reformulated optimization problem \ref{opt-gg1-r2} is convex. If moreover, Assumption~\ref{assumption:mu_hat}.\ref{assumption:mu_hat strongly convex} is satisfied, then problem \ref{opt-gg1-r2} is $\kappa_2$-strongly convex.
\end{theorem}

Theorem~\ref{theorem-GM1-logconcave} complements Theorem~\ref{theorem-ee1} by establishing the convexity of $\hat L_q(\tau)$ beyond the Gamma family to log-concave interarrival-time distributions. This result further shows that the convexity result does not rely solely on distribution-specific expressions but also follows from structural properties of the distribution.

\begin{remark} 
    The square-root transformations used in Reformulations~\ref{opt-gg1-r1} and~\ref{opt-gg1-r2} are not the only transformations that can yield convex reformulations. For any $\alpha,\beta>0$, define
    \begin{equation*}
        \hat\lambda_\alpha=\lambda^{1/\alpha},\qquad \hat\mu_\beta=\mu^\beta,\qquad \hat r_\alpha(\hat\lambda_\alpha)=r(\hat\lambda_\alpha^\alpha),\qquad \hat c_\beta(\hat\mu_\beta)=c(\hat\mu_\beta^{1/\beta}).
    \end{equation*}
    Under this joint transformation, we have
    \begin{equation*}
        L_q(\lambda/\mu) = L_q\left( \frac{\hat\lambda_\alpha^{\alpha}}{\hat\mu_\beta^{1/\beta}} \right) = \hat L_q \left( \frac{\hat\lambda_\alpha^{2\alpha}} {\hat\mu_\beta^{2/\beta}} \right), \qquad \hat L_q(\tau)=L_q(\sqrt{\tau}).
    \end{equation*}
    For the $M/M/1$ queue, the expected queue-length has the closed form $L_q(\rho)= \rho^2 / (1-\rho)$. A direct Hessian calculation shows that the transformed expected queue length is jointly convex in $(\hat\lambda_\alpha,\hat\mu_\beta)$ on the stable domain if and only if $2\alpha \ge 1+2 / \beta$. More generally, consider the following two settings:
    \begin{itemize}
        \item If $L_q$ is convex and nondecreasing, as in the setting of Theorem~\ref{theorem-gg1-r1}, then the transformed expected queue-length component remains convex in $(\hat\lambda_\alpha,\hat\mu_\beta)$ whenever $\alpha\ge 1 + 1 / \beta$. This follows from the fact that $(x,y)\mapsto x^a/y^b$ is jointly convex on the positive orthant whenever $b\ge0$ and $a\ge b+1$. In particular, this condition is satisfied by $\alpha=2,\beta=1$, which reduces to Reformulation~\ref{opt-gg1-r1}.
        
        \item If $\hat L_q$ is convex and nondecreasing, as in the settings of Theorems~\ref{theorem-ee1} and~\ref{theorem-GM1-logconcave}, then the transformed expected queue-length component remains convex in $(\hat\lambda_\alpha,\hat\mu_\beta)$ whenever $2 \alpha \ge 1+ 2 / \beta$. In particular, this condition is satisfied by $\alpha=1,\beta=2$, which reduces to Reformulation~\ref{opt-gg1-r2}. Within the family $\alpha=1,\beta\ge2$, the choice $\beta=2$ imposes the least restrictive condition on the service-cost function, because convexity of $\hat c_\beta$ with $\beta \ge 2$ implies convexity of $\hat c(\hat\mu)=c(\sqrt{\hat\mu})$ when $c$ is nondecreasing.
    \end{itemize}
    Thus, power transformations in this family can yield convex reformulations, subject to the corresponding concavity and convexity requirements on the transformed revenue and service-cost components.
\end{remark}

\subsubsection*{Counterexample: Mixture of Exponential Distributions}
Given the positive results for Gamma distributions and log-concave distributions, one may naturally wonder whether the same property holds for general phase-type distributions, i.e., whether $\hat{L}_q(\tau)$ remains convex when $S_n$ and $T_n$ are phase-type random variables. However, this is not the case, as stated in the following theorem.

\begin{theorem}\label{theorem-mix exp}
    $\hat{L}_q(\tau)$ can be non-convex on $(0,1)$ when $S_n$ and $T_n$ are general phase-type random variables.
\end{theorem}

The proof of Theorem~\ref{theorem-mix exp} follows a constructive approach. Specifically, we present a counterexample in which $S_n$ and $T_n$ are mixtures of exponential distributions, demonstrating that $\hat{L}_q(\tau)$ is not convex. By definition, it is sufficient to prove that there exists $0 < \tau_1 < \tau_2 < 1$, such that $\sum_{n=1}^\infty \hat{l}'_n(\tau_1) > \sum_{n=1}^\infty \hat{l}'_n(\tau_2)$. 

\section{Polyak-Łojasiewicz-Kurdyka (PŁK) Condition}\label{sec-PLK condition}
Building on the convex reformulations, we analyze the landscape of the original non-convex problem (\ref{opt-main problem}). Our analysis relies on a specific form of the PŁK condition (also referred to as the KŁ condition in \citealp{karimi2016linear}; see their Appendix G). We refer interested readers to \citet[Definition 1]{bento2025convergence} for the general formulation. Following \citet{bento2025convergence}, we adopt the PŁK terminology, rather than the more common PŁ or KŁ condition, to acknowledge the foundational contributions of \citet{polyak1963gradient}, \citet{lojasiewicz1963topological}, and \citet{kurdyka1998gradients}.

\begin{definition}[PŁK Condition]\label{def: PLK condition}
    Consider a convex and compact set $\mathcal{X}\subseteq\mathbb{R}^n$ and a differentiable function $f$. Denote $f^*$ as the optimal function value with $f^* \coloneqq \min_{x\in\mathcal{X}} f(x)$. The function $f$ satisfies the PŁK condition on $\cx$ if there exists $\nu>0$ such that
    \begin{equation*}
        f(x) - f^* \le \frac{1}{2\nu} \min_{g \in \partial \delta_\cx(x)} \left \| \nabla f(x) + g \right \|_2^\alpha \quad \forall x\in\cx,
    \end{equation*}
    where $\nu>0$ denotes the P{\L}K constant and $\alpha$ is a real-valued P{\L}K exponent satisfying $\alpha \ge 1$.
\end{definition}

In Definition~\ref{def: PLK condition}, $\partial \delta_\cx(x)$ denotes the subdifferential of $\delta_\cx(x)$ and is the normal cone of $\cx$ at $x$. Despite variations such as whether the domain is bounded or not, or whether the objective function is differentiable or not, an essential property of the PŁK condition is that optimization problems satisfying it exclude all suboptimal stationary or local optimal points \citep{karimi2016linear}. Hence, the PŁK condition serves as a perfect candidate to explain the empirical success of first-order methods.

\begin{theorem}\label{theorem-PLK}
    Suppose Assumption~\ref{assumption-general model} holds. Consider problem (\ref{opt-main problem}) in one of the following settings:
    \begin{henumerate}
        \item If Assumption~\ref{assumption:lambda_hat}.\ref{assumption:lambda_hat convex} holds, then (\ref{opt-main problem}) satisfies the PŁK condition with exponent $\alpha = 1$. Moreover, if Assumption~\ref{assumption:lambda_hat}.\ref{assumption:lambda_hat strongly convex} is satisfied, then (\ref{opt-main problem}) satisfies the PŁK condition with exponent $\alpha = 2$.
        \item Suppose that $T_n\sim\Gamma(k, k)$ and $S_n\sim\Gamma(m, m)$ with $k, m \ge 1$. If Assumption~\ref{assumption:mu_hat}.\ref{assumption:mu_hat convex} holds, then (\ref{opt-main problem}) satisfies the PŁK condition with exponent $\alpha = 1$. Moreover, if Assumption~\ref{assumption:mu_hat}.\ref{assumption:mu_hat strongly convex} is satisfied, then (\ref{opt-main problem}) satisfies the PŁK condition with exponent $\alpha = 2$.
    \end{henumerate}
\end{theorem}

The intuition behind Theorem~\ref{theorem-PLK} is as follows. The original problem \eqref{opt-main problem} can be viewed as the reformulation composed with a transformation map of decision variables. For simplicity, let us express \eqref{opt-main problem} by:
\begin{equation}
\label{function composition}
    \begin{aligned}
        \min_{x \in \cx}\ F(x)=H\left(\varphi(x)\right),
    \end{aligned}
\end{equation}
where $F$ and $\cx$ are the objective function and feasible region of \eqref{opt-main problem}, respectively, $x$ denotes the decision variables $(\lambda,\mu)$, $\varphi(x)$ is the transformation map used to construct the reformulation, and $H$ is the objective function of the reformulation. \citet[Theorem 3.2]{li2018calculus} and \citet[Proposition 2]{fatkhullin2023stochastic} prove that the P{\L}K condition can be established through such a composition when the transformation map is uniformly inverse-Lipschitz on the feasible region. Specifically, if the reformulated feasible region $\cu=\varphi(\cx)$ is closed and convex, and there exists $\sigma_\varphi>0$ such that
\begin{equation*}
    \|\varphi(x)-\varphi(y)\|_2 \ge \sigma_\varphi \|x-y\|_2, \qquad \forall x,y\in\cx,
\end{equation*}
then convexity of $H$ on $\cu$ implies that $F=H\circ\varphi$ satisfies the P{\L}K condition with exponent $\alpha=1$ on $\cx$. Moreover, if $H$ is strongly convex on $\cu$, then $F$ satisfies the P{\L}K condition with exponent $\alpha=2$.

In our case, we have established the convexity and strong convexity results for reformulations \ref{opt-gg1-r1} and \ref{opt-gg1-r2} in Theorem~\ref{theorem-gg1-r1} and Theorem~\ref{theorem-ee1}, respectively. In the proof of Theorem~\ref{theorem-PLK}, we verify the required inverse-Lipschitz condition directly for the two transformation maps $\varphi_1(\lambda,\mu)=(\sqrt{\lambda},\mu)$ and $\varphi_2(\lambda,\mu)=(\lambda,\mu^2)$. 

\begin{remark}
    For the Jackson networks, the $M/GI/1$ queue, the $GI/GI/s$ approximation, and the $GI/M/1$ queue with log-concave interarrival times discussed in Section~\ref{section: reformulation II}, the proof extends directly to establish the P\L K condition with exponent $\alpha=1$ under Assumption~\ref{assumption:mu_hat}.\ref{assumption:mu_hat convex}, and with exponent $\alpha=2$ under Assumption~\ref{assumption:mu_hat}.\ref{assumption:mu_hat strongly convex}.
\end{remark}

Theorem~\ref{theorem-PLK} shows that the original optimization problem satisfies the P{\L}K condition, providing a strong characterization of the non-convex landscape. This insight also extends to formulations in which the decision variables are price $p$ and service rate $\mu$, rather than arrival rate $\lambda$ and service rate $\mu$. Specifically, with the demand model $\lambda=\lambda(p)$, where $\lambda(\cdot)$ is continuously differentiable and $\inf_{p\in\mathcal P} |\lambda'(p)| > 0$ on a compact price interval $\mathcal{P}$, the mapping $(p,\mu)\mapsto (\lambda(p),\mu)$ is uniformly inverse-Lipschitz. Hence, by the same argument used in the proof of Theorem~\ref{theorem-PLK}, the P{\L}K condition is preserved under this change of variables, up to a modification of the constant. Therefore, the formulation with decision variables $(p,\mu)$ also satisfies a P{\L}K condition. This shows that our theory directly applies to the practically relevant setting in which the decision maker controls price rather than arrival rate, and thus provides a theoretical explanation for the success of first-order methods in the $(p,\mu)$ formulations studied by \citet{chen2024online} and \citet{hu2024multi}.

\section{Technical Proofs}\label{sec-proof}
\subsection{Proof of Theorem~\ref{theorem-ee1}}
\begin{proof}{Proof of Theorem~\ref{theorem-ee1}.}
    For the proof, it is convenient to rescale the unit-mean Gamma variables. Let $\widetilde T_n \coloneqq kT_n\sim\Gamma(k,1)$ and $\widetilde S_n \coloneqq mS_n\sim\Gamma(m,1)$, and let $\widetilde L_q$ denote the corresponding transformed queue-length function. Then, $\hat L_q(\tau)=\frac{1}{k}\widetilde L_q (\tau k^2 / m^2 )$. Therefore, it suffices to prove that $\widetilde L_q$ is convex on $(0,k^2/m^2)$. For notational simplicity, we write $\hat L_q$ for $\widetilde L_q$ and suppress the tildes on $S_n$ and $T_n$ throughout the proof.

    We first establish a monotonicity property needed later. For $n\ge1$, define
    \begin{equation*}
        s_n \coloneqq \frac{\Gamma((k+m)n)}{\Gamma(mn)\Gamma(kn)} \left[\frac{m^mk^k}{(m+k)^{m+k}}\right]^n.
    \end{equation*}
    We claim that $s_n$ is non-decreasing in $n$. To see this, define, for $x>0$,
    \begin{equation*}
        h(x) \coloneqq \ln\Gamma((k+m)x)-\ln\Gamma(kx)-\ln\Gamma(mx) -x\left[k\ln\left(\frac{k+m}{k}\right)+m\ln\left(\frac{k+m}{m}\right)\right].
    \end{equation*}
    Therefore, $s_n=e^{h(n)}$. Let $\psi(x)\coloneqq \Gamma'(x)/\Gamma(x)$ be the digamma function. Using its integral representation \citep[Equation (6.3.22) on Page 259]{abramowitz1965handbook}, we obtain
    \begin{equation*}
        \begin{aligned}
            h'(x) &= k\int_0^\infty \frac{e^{-kxt}-e^{-(k+m)xt}}{1-e^{-t}}dt +m\int_0^\infty \frac{e^{-mxt}-e^{-(k+m)xt}}{1-e^{-t}}dt \\
            &\quad -k\ln\left(\frac{k+m}{k}\right)-m\ln\left(\frac{k+m}{m}\right).
        \end{aligned}
    \end{equation*}
    Since $1-e^{-t}\le t$ and the numerators in the two integrals are non-negative, Frullani's integral \citep[Equation (3.434.2) on Page 363]{gradshteyn2014table} gives
    \begin{equation*}
        \begin{aligned}
            h'(x) &\ge k\int_0^\infty \frac{e^{-kxt}-e^{-(k+m)xt}}{t}dt +m\int_0^\infty \frac{e^{-mxt}-e^{-(k+m)xt}}{t}dt \\
            &\quad -k\ln\left(\frac{k+m}{k}\right)-m\ln\left(\frac{k+m}{m}\right)=0.
        \end{aligned}
    \end{equation*}
    Thus, $s_n$ is non-decreasing in $n$.

    Now define
    \begin{equation*}
        U_n\coloneqq \sum_{i=1}^n(S_i+T_i), \qquad V_n\coloneqq \frac{\sum_{i=1}^n S_i}{\sum_{i=1}^n(S_i+T_i)}.
    \end{equation*}
    By the standard beta-gamma decomposition, $U_n\sim\Gamma((k+m)n,1)$, $V_n\sim\text{Beta}(mn,kn)$, and $U_n$ is independent of $V_n$. Let $p_n$ denote the density of $V_n$, i.e.,
    \begin{equation*}
        p_n(v)=\frac{\Gamma((k+m)n)}{\Gamma(mn)\Gamma(kn)}v^{mn-1}(1-v)^{kn-1}.
    \end{equation*}
    Hence,
    \begin{equation*}
        \begin{aligned}
            \hat{l}_n(\tau) = \frac{1}{n} \E \left[ \left(\sum_{i=1}^n \sqrt{\tau} S_i - \sum_{i=1}^n T_i \right)^+ \right] = \frac{1}{n}\E \left[ U_n \left( (\sqrt{\tau} + 1) V_n - 1 \right)^+ \right] = (k+m)\E\left[ \left( (\sqrt{\tau} + 1) V_n - 1 \right)^+ \right].
        \end{aligned}
    \end{equation*}
    It follows that
    \begin{equation}\label{eq-hat-L-gamma-compact}
        \hat L_q(\tau) =(k+m)\sum_{n=1}^\infty\int_a^1\left[(1+\sqrt{\tau})v-1\right]p_n(v)dv, \qquad a\coloneqq\frac{1}{1+\sqrt{\tau}}.
    \end{equation}
    Since $\tau\in(0,k^2/m^2)$, we have $a\in(m/(m+k),1)$.

    The termwise differentiations below are justified by uniform convergence on compact subintervals of $(0,k^2/m^2)$. This follows because, on such intervals, the differentiated terms are bounded by a summable geometrically decaying sequence. Differentiating \eqref{eq-hat-L-gamma-compact} twice and using $a=(1+\sqrt{\tau})^{-1}$, we obtain
    \begin{equation}\label{eq-gamma-second-derivative-D}
        \hat L_q''(\tau) =\frac{(k+m)a^3}{4(1-a)^3}D(a), \qquad D(a)\coloneqq
        \sum_{n=1}^\infty\left[a^2(1-a)p_n(a)-\int_a^1 vp_n(v)dv\right].
    \end{equation}
    Hence, it remains to prove that $D(a)\ge0$ for all $a\in(m/(m+k),1)$.

    Observe that $D(1)=0$. Thus, it is sufficient to prove that $D'(a)\le0$ on $(m/(m+k),1)$. Since
    \begin{equation*}
        p_n'(a)=\frac{mn-(k+m)an-(1-2a)}{a(1-a)}p_n(a),
    \end{equation*}
    direct differentiation gives
    \begin{equation}\label{eq-gamma-D-prime-compact}
        D'(a) =(2a-a^2)\sum_{n=1}^\infty p_n(a) -a\left[(k+m)a-m\right]\sum_{n=1}^\infty np_n(a).
    \end{equation}
    Define
    \begin{equation*}
        r(a) \coloneqq \frac{(k+m)^{k+m}}{k^km^m}a^m(1-a)^k.
    \end{equation*}
    For $a\in(m/(m+k),1)$, we have $0<r(a)<1$, since $a^m(1-a)^k$ attains its maximum at $m/(m+k)$. Moreover,
    \begin{equation*}
        p_n(a)= \frac{(k+m)^{k+m}}{k^km^m}a^{m-1}(1-a)^{k-1}s_n\left[r(a)\right]^{n-1}.
    \end{equation*}
    Since $s_n$ is non-decreasing in $n$ and $0<r(a)<1$,
    \begin{equation}\label{eq-gamma-series-lower-bound}
        \sum_{n=1}^\infty n p_n(a) \ge \frac{1}{1-r(a)}\sum_{n=1}^\infty p_n(a).
    \end{equation}
    Indeed, after removing the common positive factor in $p_n(a)$, this follows from
    \begin{equation*}
        \sum_{n=1}^\infty ns_n r^{n-1} =\sum_{q=1}^\infty\sum_{n=q}^\infty s_n r^{n-1} \ge \sum_{q=1}^\infty\sum_{n=q}^\infty s_q r^{n-1} =\frac{1}{1-r}\sum_{q=1}^\infty s_q r^{q-1},
    \end{equation*}
    where $r=r(a)$.

    Combining \eqref{eq-gamma-D-prime-compact} and \eqref{eq-gamma-series-lower-bound}, we obtain
    \begin{equation}\label{eq-gamma-D-prime-upper-bound}
        D'(a) \le a\left(\sum_{n=1}^\infty p_n(a)\right) \left[2-a-\frac{(k+m)a-m}{1-r(a)}\right].
    \end{equation}
    Therefore, it remains only to prove the scalar inequality
    \begin{equation}\label{eq-gamma-kernel-bound-compact}
        (k+m)a-m\ge (2-a)(1-r(a)), \qquad a\in\left(\frac{m}{m+k},1\right).
    \end{equation}
    If $a\ge (m+2)/(m+k+1)$, then $(k+m)a-m\ge 2-a\ge (2-a)(1-r(a))$, so \eqref{eq-gamma-kernel-bound-compact} holds. It remains to consider $m/(m+k)<a<(m+2)/(m+k+1)$. Define
    \begin{equation*}
        \eta\coloneqq\frac{(k+m)a-m}{2-a}\in(0,1), \qquad a=\frac{m+2\eta}{k+m+\eta}.
    \end{equation*}
    Then \eqref{eq-gamma-kernel-bound-compact} is equivalent to $r(a)\ge1-\eta$. Substituting the above expression for $a$ into $r(a)$ gives
    \begin{equation*}
        r(a)= \left(\frac{k+m}{k+m+\eta}\right)^{k+m} \left(1+\frac{2\eta}{m}\right)^m \left(1-\frac{\eta}{k}\right)^k.
    \end{equation*}
    Let $s\coloneqq k+m$. Since $0<2\eta/m<2$ and $\ln(1+z)\ge z/2$ for $0\le z\le2$, we have $(1+2\eta/m)^m\ge e^{\eta}$. Also, $(1+\eta/s)^s\le e^{\eta}$. Therefore, $(1+ 2\eta / m )^m \ge (1+ \eta / s)^s$. In addition, Bernoulli's inequality gives $(1-\eta/k)^k\ge1-\eta$. Hence,
    \begin{equation*}
        r(a) \ge \left(\frac{s}{s+\eta}\right)^s \left(1+\frac{\eta}{s}\right)^s (1-\eta) =1-\eta.
    \end{equation*}
    This proves \eqref{eq-gamma-kernel-bound-compact}. By \eqref{eq-gamma-D-prime-upper-bound}, $D'(a)\le0$ on $(m/(m+k),1)$. Since $D(1)=0$, we have $D(a)\ge0$ on $(m/(m+k),1)$. Therefore, \eqref{eq-gamma-second-derivative-D} implies that $\hat L_q$ is convex on $(0,k^2/m^2)$.

    The preceding argument proves that $\widetilde L_q$ is convex and nondecreasing on $(0,k^2/m^2)$. Therefore, for the original unit-mean Gamma variables, $\hat L_q(\tau)=\frac{1}{k}\widetilde L_q (\tau k^2 / m^2 )$ is convex and nondecreasing on $(0,1)$. Under Assumption~\ref{assumption-general model}, $\pi(\lambda,\hat\mu)\in(0,1)$ for all feasible $(\lambda,\hat\mu)$. Since $\hat L_q$ is convex and nondecreasing on $(0,1)$, and $\pi(\lambda,\hat\mu)$ is convex, the composition rule implies that $\hat L_q(\pi(\lambda,\hat\mu))$ is jointly convex in $(\lambda,\hat\mu)$. Under Assumption~\ref{assumption:mu_hat}.\ref{assumption:mu_hat convex}, the function $\hat c(\hat\mu)-r(\lambda)$ is jointly convex in $(\lambda,\hat\mu)$. Therefore, problem \ref{opt-gg1-r2} is convex. Moreover, under Assumption~\ref{assumption:mu_hat}.\ref{assumption:mu_hat strongly convex}, the function $\hat c(\hat\mu)-r(\lambda)$ is jointly $\kappa_2$-strongly convex in $(\lambda,\hat\mu)$. Because the sum of a strongly convex function and a convex function remains strongly convex with the same modulus, problem \ref{opt-gg1-r2} is $\kappa_2$-strongly convex. \Halmos
\end{proof}

\subsection{Proof of Theorem~\ref{theorem-GM1-logconcave}}
\begin{proof}{Proof of Theorem~\ref{theorem-GM1-logconcave}.}
    Let $\phi(u)\coloneqq \E [ e^{-uT_1} ]$ with $u > 0$ be the Laplace transform of the interarrival-time distribution. For the $GI/M/1$ queue with traffic intensity $\rho\in(0,1)$, let $\sigma\in(0,1)$ be the non-trivial solution to $\sigma=\phi ( (1-\sigma)/\rho )$. The stationary waiting time in queue satisfies $\P(W_\infty>t)=\sigma e^{-\mu(1-\sigma)t}$ for $t \ge 0$. Hence, by Little's law, $L_q(\rho)=\lambda\E[W_\infty] = \rho\sigma / (1-\sigma)$. Define $u \coloneqq (1-\sigma) / \rho$. Then,
    \begin{equation*}
        \sigma=\phi(u), \qquad \rho(u)=\frac{1-\phi(u)}{u}, \qquad \hat L_q(\tau(u))=\frac{\phi(u)}{u}=\frac{1}{u}-\rho(u), \qquad \tau(u)=\rho(u)^2.
    \end{equation*}
    Let $\bar F(t)\coloneqq \P(T_1>t)$. Since $\E[T_1]=1$, integration by parts gives
    \begin{equation}\label{eq-GM1-rho-integral}
        \rho(u)=\frac{1-\phi(u)}{u} = \int_0^\infty e^{-ut}\bar F(t)dt.
    \end{equation}
    Hence,
    \begin{equation*}
        \rho'(u)=-\int_0^\infty te^{-ut}\bar F(t)dt<0,\qquad u>0.
    \end{equation*}
    Moreover, $\lim_{u\to 0^+}\rho(u)=\E[T_1]=1$ and $\lim_{u\to\infty} \rho(u) = 0$. Thus, $u\mapsto\rho(u)$ is a bijection from $(0,\infty)$ to $(0,1)$, and $u\mapsto\tau(u)$ is a bijection from $(0,\infty)$ to $(0,1)$.

    Define $C(u) \coloneqq - \rho'(u) / \rho(u)$. By direct differentiation of the parametric representation $\hat L_q(\tau(u))=1/u-\rho(u)$ and $\tau(u)=\rho(u)^2$, we obtain
    \begin{equation}\label{eq-GM1-second-derivative-rho}
        \hat L_q''(\tau(u))
        =
        \frac{2C(u)-2uC(u)^2+uC'(u)+u^3\rho(u)C(u)^3}
        {4u^3\rho(u)^4C(u)^3}.
    \end{equation}
    Since the denominator in \eqref{eq-GM1-second-derivative-rho} is positive, it remains to prove that the numerator is positive.

    We next establish two bounds for $C(u)$. Let $f$ be the density of $T_1$, extended by zero outside its support. Since $f$ is log-concave on an interval, and log-concavity is preserved under integration, the survival function $\bar F(t)$ is log-concave on $[0,\infty)$. Therefore, for each fixed $u>0$, the density
    \begin{equation*}
        g_u(t)\coloneqq \frac{e^{-ut}\bar F(t)}{\rho(u)}, \qquad t\ge0,
    \end{equation*}
    is log-concave. Let $Z_u$ denote a random variable with density $g_u$. By \eqref{eq-GM1-rho-integral},
    \begin{equation*}
        \E[Z_u]=C(u), \qquad C'(u)=-\mathrm{Var}(Z_u).
    \end{equation*}

    Let $\bar G_u$ and $r_u$ denote the survival function and hazard rate of $Z_u$, respectively. Since $g_u$ is log-concave, $\bar G_u$ is log-concave, and hence $r_u$ is increasing. Moreover, $\bar F(0)=1$, so $r_u(0)=g_u(0)=1 / \rho(u)$. Hence,
    \begin{equation*}
        \bar G_u(t) = \exp\left\{-\int_0^t r_u(s)ds\right\} \le e^{-t/\rho(u)},  \qquad \forall \ t\ge0.
    \end{equation*}
    Integrating both sides gives
    \begin{equation}\label{eq-GM1-C-leq-rho}
        C(u)=\E[Z_u]=\int_0^\infty \bar G_u(t)dt\le \rho(u).
    \end{equation}

    We also need a lower bound on $C'(u)$. Since $r_u$ is increasing, for any $s,t\ge0$,
    \begin{equation*}
        \frac{\bar G_u(t+s)}{\bar G_u(t)} = \exp\left\{-\int_t^{t+s}r_u(v)dv\right\} \le \exp\left\{-\int_0^s r_u(v)dv\right\} = \bar G_u(s).
    \end{equation*}
    Thus,
    \begin{equation*}
        \int_t^\infty \bar G_u(s)ds = \bar G_u(t)\int_0^\infty \frac{\bar G_u(t+s)}{\bar G_u(t)}ds \le C(u)\bar G_u(t).
    \end{equation*}
    It follows that
    \begin{equation*}
        \begin{aligned}
            \E[Z_u^2] =
            2\int_0^\infty t\bar G_u(t)dt =
            2\int_0^\infty\int_t^\infty \bar G_u(s)dsdt \le
            2C(u)\int_0^\infty \bar G_u(t)dt = 2C(u)^2.
        \end{aligned}
    \end{equation*}
    Hence,
    \begin{equation}\label{eq-GM1-Cprime-lower-rho}
        C'(u)=-\mathrm{Var}(Z_u)\ge -C(u)^2.
    \end{equation}

    Combining \eqref{eq-GM1-second-derivative-rho} with \eqref{eq-GM1-Cprime-lower-rho}, it is sufficient to show that $2C(u)-3uC(u)^2+u^3\rho(u)C(u)^3>0$. Set
    \begin{equation*}
        x\coloneqq uC(u), \qquad y\coloneqq u\rho(u)=1-\phi(u).
    \end{equation*}
    Then, $x>0$, $0<y<1$, and $x\le y$ by \eqref{eq-GM1-C-leq-rho}. Moreover,
    \begin{equation*}
        2C(u)-3uC(u)^2+u^3\rho(u)C(u)^3 = C(u)\left(2-3x+yx^2\right).
    \end{equation*}
    For any fixed $y\in(0,1)$, the function $2-3x+yx^2$ is decreasing in $x$ on $(0,y]$, because $2yx-3\le2y^2-3<0$. Hence, $2-3x+yx^2\ge y^3-3y+2=(1-y)^2(y+2)>0$. Therefore, $\hat L_q''(\tau(u))>0$ for every $u>0$. Since $u\mapsto\tau(u)$ is a bijection from $(0,\infty)$ to $(0,1)$, $\hat L_q(\tau)$ is convex on $(0,1)$.

    Moreover, $\hat L_q$ is nondecreasing on $(0,1)$ because $L_q(\rho)$ is nondecreasing in $\rho$. Under Assumption~\ref{assumption-general model}, $\pi(\lambda,\hat\mu)\in(0,1)$ for all feasible $(\lambda,\hat\mu)$. Since $\hat L_q$ is convex and nondecreasing on $(0,1)$, and $\pi(\lambda,\hat\mu)$ is convex, the composition rule implies that $\hat L_q(\pi(\lambda,\hat\mu))$ is jointly convex in $(\lambda,\hat\mu)$. Under Assumption~\ref{assumption:mu_hat}.\ref{assumption:mu_hat convex}, the function $\hat c(\hat\mu)-r(\lambda)$ is jointly convex in $(\lambda,\hat\mu)$. Therefore, problem \ref{opt-gg1-r2} is convex. Moreover, under Assumption~\ref{assumption:mu_hat}.\ref{assumption:mu_hat strongly convex}, the function $\hat c(\hat\mu)-r(\lambda)$ is jointly $\kappa_2$-strongly convex in $(\lambda,\hat\mu)$. Because the sum of a strongly convex function and a convex function remains strongly convex, problem \ref{opt-gg1-r2} is $\kappa_2$-strongly convex. \Halmos
\end{proof}

\subsection{Proof of Theorem~\ref{theorem-mix exp}}
\begin{proof}{Proof of Theorem~\ref{theorem-mix exp}.}
    Consider a setting when $S_n$ and $T_n$ are mixtures of exponential random variables:
    \begin{equation}\label{rv-counterexample}
        S_n, T_n \overset{\text{i.i.d.}}{\sim} \left\{
        \begin{aligned}
        & \text{Exp}(\lambda_1 = 0.1), &\text{w.p.} \ \frac{1}{2};\\
        & \text{Exp}(\lambda_2 = 10), &\text{w.p.} \ \frac{1}{2}.\\
        \end{aligned}
        \right.
    \end{equation}
    To enforce $\E[S_n] = \E[T_n] = 1$, we should use $S_n / \E[S_n]$ and $T_n / \E[T_n]$. In our example, $\E[S_n] = \E[T_n]$, so both sequences are rescaled by the same positive constant, which does not affect the non-convexity result. For simplicity, we keep using $S_n$ and $T_n$ as in \eqref{rv-counterexample} in the following analysis.
    
    Let us define $X_n \coloneqq \sum_{i=1}^n S_i$ and $Y_n \coloneqq \sum_{i=1}^n T_i$. By definition,
    \begin{equation*}
        \begin{aligned}
            \hat{l}_n(\tau) = \frac{1}{n} \E \left[ \left( \sqrt{\tau}X_n - Y_n \right)^+ \right] = \frac{1}{n} \int_0^\infty \left( \int_0^{\sqrt{\tau}x} (\sqrt{\tau}x - y) p_{X_n}(x)p_{Y_n}(y) dy\right) dx.
        \end{aligned}
    \end{equation*}
    Applying the Leibniz rule:
    \begin{equation*}
        \begin{aligned}
            \hat{l}'_n(\tau) = \frac{1}{2n\sqrt{\tau}} \int_0^\infty \left( \int_0^{\sqrt{\tau}x} x p_{X_n}(x)p_{Y_n}(y) dy\right) dx.
        \end{aligned}
    \end{equation*}
    Set $\tau_1 = 10^{-6}$ and $\tau_2 = 0.01$. It is sufficient to prove that $\sum_{n=1}^\infty \hat{l}'_n(\tau_1) > \sum_{n=1}^\infty \hat{l}'_n(\tau_2)$.

    \textbf{Step 1: Lower bound of $\sum_{n=1}^\infty \hat{l}'_n(\tau)$}: Since $\hat{l}'_n(\tau)$ is non-negative, we have:
    \begin{equation*}
        \sum_{n=1}^\infty \hat{l}'_n(\tau) \ge \hat{l}'_1(\tau),\quad \forall \tau \in (0,1).
    \end{equation*}
    When $n=1$, $X_1 = S_1$ and $Y_1 = T_1$ with PDFs:
    \begin{equation*}
        \begin{aligned}
            p_{S_1}(s) = \frac{\lambda_1}{2} e^{-\lambda_1 s} + \frac{\lambda_2}{2}e^{-\lambda_2 s}, \quad p_{T_1}(t) = \frac{\lambda_1}{2} e^{-\lambda_1 t} + \frac{\lambda_2}{2}e^{-\lambda_2 t}.
        \end{aligned}
    \end{equation*}
    By definition, we have:
    \begin{equation}\label{eq-mix Erlang}
        \begin{aligned}
            \hat{l}'_1(\tau) &= \frac{1}{2\sqrt{\tau}} \int_0^\infty s p_{S_1}(s) \left( \int_0^{\sqrt{\tau}s} p_{T_1}(t) dt\right) ds\\
            &= \frac{1}{8\sqrt{\tau}} \sum_{i,j=1}^2 \lambda_i \int_0^\infty s e^{-\lambda_i s} \left(1 - e^{-\lambda_j\sqrt{\tau} s} \right) ds\\
            &= \frac{1}{8\sqrt{\tau}} \sum_{i,j=1}^2 \left[ \frac{1}{\lambda_i} - \frac{\lambda_i}{(\lambda_i + \lambda_j \sqrt{\tau})^2} \right].
        \end{aligned}
    \end{equation}
    Plugging $\lambda_1 = 0.1$, $\lambda_2 = 10$, and $\tau_1 = 10^{-6}$ into the formulation, we have:
    \begin{equation*}
        \sum_{n=1}^\infty \hat{l}'_n(\tau_1) \ge \hat{l}'_1(\tau_1) > 219.
    \end{equation*}

    \textbf{Step 2: Upper bound of $\sum_{n=1}^\infty \hat{l}'_n(\tau)$}: We split the summation into two parts: $n=1$ and $n\ge2$. When $n = 1$, from \eqref{eq-mix Erlang}, we already have
    \begin{equation*}
        \hat{l}'_1(\tau) = \frac{1}{8\sqrt{\tau}} \sum_{i=1}^2\sum_{j=1}^2 \left[ \frac{1}{\lambda_i} - \frac{\lambda_i}{(\lambda_i + \lambda_j \sqrt{\tau})^2} \right].
    \end{equation*}
    Plugging $\lambda_1 = 0.1$, $\lambda_2 = 10$, and $\tau_2 = 0.01$ into the above formulation, we have $\hat{l}'_1(\tau_2) < 15$. When $n \ge 2$,
    \begin{equation*}
        \begin{aligned}
            \hat{l}'_n(\tau) = \frac{1}{2n\sqrt{\tau}} \int_0^\infty x p_{X_n}(x) \left( \int_0^{\sqrt{\tau}x} p_{Y_n}(y) dy\right) dx = \frac{1}{2n\sqrt{\tau}} \E \left[ X_n \mathbf{1}_{(Y_n \le \sqrt{\tau} X_n)} \right].
        \end{aligned}
    \end{equation*}
    For any $\alpha > 0$, we have $\mathbf{1}_{(Y_n \le \sqrt{\tau} X_n)} \le e^{\alpha(\sqrt{\tau} X_n - Y_n)}$. Thus, we have:
    \begin{equation*}
        \begin{aligned}
            \hat{l}'_n(\tau) &\le \frac{1}{2n\sqrt{\tau}} \E \left[ X_n e^{\alpha(\sqrt{\tau} X_n - Y_n)} \right]\\
            &= \frac{1}{2n\sqrt{\tau}} \E \left[ \left( \sum_{i=1}^n S_i \right) \prod_{i=1}^n e^{\alpha(\sqrt{\tau} S_i - T_i)} \right]\\
            &\le \frac{1}{2n\sqrt{\tau}} \E \left[ \sum_{i=1}^n \left( S_i e^{\alpha\sqrt{\tau}S_i} \prod_{j\neq i} e^{\alpha(\sqrt{\tau} S_j - T_j)} \right) \right],
        \end{aligned}
    \end{equation*}
    where the last inequality holds because $e^{-\alpha T_i} \le 1$. Since $\{S_i\}_{i=1}^n$ are i.i.d. and independent of $\{T_i\}_{i=1}^n$, we have
    \begin{equation*}
        \begin{aligned}
            \hat{l}'_n(\tau)\leq \frac{1}{2\sqrt{\tau}} \E \left[ S_1 e^{\alpha \sqrt{\tau} S_1} \right] \left( M_S(\alpha\sqrt{\tau}) M_T(-\alpha) \right)^{n-1},
        \end{aligned}
    \end{equation*}
    where $M_S$ and $M_T$ are the moment-generating functions of $S_i$ and $T_i$. Set $\alpha = 0.2$ and $\tau_2 = 0.01$, we can calculate each term one by one:
    \begin{equation*}
        \begin{aligned}
            \E \left[ S_1 e^{\alpha \sqrt{\tau_2} S_1} \right] &= \frac{1}{2}\sum_{i=1}^2 \int_0^\infty s e^{\alpha \sqrt{\tau_2} s} \lambda_i e^{-\lambda_i s} ds = \frac{1}{2}\sum_{i=1}^2 \frac{\lambda_i}{(\lambda_i - \alpha \sqrt{\tau_2})^2} < 8,\\
            M_S(\alpha \sqrt{\tau_2}) &= \frac{1}{2}\sum_{i=1}^2 \frac{\lambda_i}{\lambda_i - \alpha \sqrt{\tau_2}} < 1.13,\\
            M_T(-\alpha) &= \frac{1}{2}\sum_{i=1}^2 \frac{\lambda_i}{\lambda_i + \alpha} < 0.66.
        \end{aligned}
    \end{equation*}
    Combining all the results, we have $\hat{l}'_n(\tau_2) < 40 \times (0.75)^{n-1}$. Therefore, we have
    \begin{equation*}
        \begin{aligned}
            \sum_{n=1}^\infty \hat{l}'_n(\tau_2) = \hat{l}'_1(\tau_2) + \sum_{n=2}^\infty \hat{l}'_n(\tau_2) < 15 + 40 \sum_{n=2}^\infty (0.75)^{n-1} = 135.
        \end{aligned}
    \end{equation*}

    As a conclusion, we have $\sum_{n=1}^\infty \hat{l}'_n(\tau_1) > 219 > 135 > \sum_{n=1}^\infty \hat{l}'_n(\tau_2)$ with $\tau_1 = 10^{-6} < 0.01 = \tau_2$, which violates the convexity of $\hat{L}_q(\tau)$. \Halmos
\end{proof}

\subsection{Proof of Theorem~\ref{theorem-PLK}}

\begin{lemma}[{\citealp{li2018calculus, fatkhullin2023stochastic}}]\label{lemma-PLK composition}
    Consider the composition problem \eqref{function composition}. Let $\cx$ be closed and convex, and let $\cu \coloneqq \varphi(\cx)$ denote the feasible region of the reformulated problem. Suppose that $\cu$ is also closed and convex. Assume further that there exists a constant $\sigma_\varphi > 0$ such that
    \begin{equation*}
        \|\varphi(x)-\varphi(y)\|_2 \ge \sigma_\varphi \|x-y\|_2, \qquad \forall x,y \in \cx.
    \end{equation*}
    Then, the following statements hold.
    \begin{henumerate}
        \item If $H$ is convex on $\cu$ and $\cu$ is bounded with diameter $D_{\cu}$, then $F$ satisfies the PŁK condition on $\cx$ with exponent $\alpha = 1$ and constant $\sigma_\varphi/(2D_{\cu})$.
        \item If $H$ is $\alpha_H$-strongly convex on $\cu$, then $F$ satisfies the PŁK condition on $\cx$ with exponent $\alpha = 2$ and constant $\alpha_H \sigma_\varphi^2$.
    \end{henumerate}
\end{lemma}

\begin{proof}{Proof of Theorem \ref{theorem-PLK}}
    First consider reformulation \ref{opt-gg1-r1}, which uses the change of variables $\hat\lambda=\sqrt{\lambda}$. Set $\varphi_1(\lambda,\mu)=(\sqrt{\lambda},\mu)$. Then, $\cu_1=\varphi_1(\cx)=\{(\hat\lambda,\mu): \sqrt{\underline\lambda}\le \hat\lambda\le \sqrt{\bar\lambda},\ \underline\mu\le \mu\le \bar\mu\}$, which is closed and convex. Moreover, $\cu_1$ is bounded with diameter $D_{\cu_1}=\sqrt{(\sqrt{\bar\lambda}-\sqrt{\underline\lambda})^2+(\bar\mu-\underline\mu)^2}$. For any $(\lambda_1,\mu_1),(\lambda_2,\mu_2)\in\cx$,
    \begin{equation*}
        \|\varphi_1(\lambda_1,\mu_1)-\varphi_1(\lambda_2,\mu_2)\|_2^2
    =(\sqrt{\lambda_1}-\sqrt{\lambda_2})^2+(\mu_1-\mu_2)^2.
    \end{equation*}
    Since
    \begin{equation*}
        |\sqrt{\lambda_1}-\sqrt{\lambda_2}| =\frac{|\lambda_1-\lambda_2|}{\sqrt{\lambda_1}+\sqrt{\lambda_2}} \ge \frac{|\lambda_1-\lambda_2|}{2\sqrt{\bar\lambda}},
    \end{equation*}
    we obtain
    \begin{equation*}
        \|\varphi_1(\lambda_1,\mu_1)-\varphi_1(\lambda_2,\mu_2)\|_2^2 \ge \frac{1}{4\bar\lambda}(\lambda_1-\lambda_2)^2+(\mu_1-\mu_2)^2 \ge \min\left\{\frac{1}{4\bar\lambda},1\right\}\|(\lambda_1,\mu_1)-(\lambda_2,\mu_2)\|_2^2.
    \end{equation*}
    By Theorem~\ref{theorem-gg1-r1}, the objective of reformulation \ref{opt-gg1-r1} is convex under Assumption~\ref{assumption:lambda_hat}.\ref{assumption:lambda_hat convex}, and $\kappa_1$-strongly convex under Assumption~\ref{assumption:lambda_hat}.\ref{assumption:lambda_hat strongly convex}. Therefore, Lemma~\ref{lemma-PLK composition} implies that the original problem \eqref{opt-main problem} satisfies the P\L K condition with exponent $\alpha=1$ and constant
    \begin{equation*}
        \frac{\min\{1/(2\sqrt{\bar\lambda}), 1\}} {2\sqrt{(\sqrt{\bar\lambda}-\sqrt{\underline\lambda})^2+(\bar\mu-\underline\mu)^2}}
    \end{equation*}
    under Assumption~\ref{assumption:lambda_hat}.\ref{assumption:lambda_hat convex}, and satisfies the P\L K condition with exponent $\alpha=2$ and constant $\kappa_1\min \{1 / (4\bar\lambda),1 \}$ under Assumption~\ref{assumption:lambda_hat}.\ref{assumption:lambda_hat strongly convex}.

    Next consider reformulation \ref{opt-gg1-r2}, which uses the change of variables $\hat\mu=\mu^2$. Set $\varphi_2(\lambda,\mu)=(\lambda,\mu^2)$. Then, $\cu_2=\varphi_2(\cx)=\{(\lambda,\hat\mu): \underline\lambda\le \lambda\le \bar\lambda,\ \underline\mu^2\le \hat\mu\le \bar\mu^2\}$, which is closed and convex. Moreover, $\cu_2$ is bounded with diameter $D_{\cu_2}=\sqrt{(\bar\lambda-\underline\lambda)^2+(\bar\mu^2-\underline\mu^2)^2}$. For any $(\lambda_1,\mu_1),(\lambda_2,\mu_2)\in\cx$,
    \begin{equation*}
        \|\varphi_2(\lambda_1,\mu_1)-\varphi_2(\lambda_2,\mu_2)\|_2^2 =(\lambda_1-\lambda_2)^2+(\mu_1^2-\mu_2^2)^2.
    \end{equation*}
    Because
    \begin{equation*}
        |\mu_1^2-\mu_2^2| =|\mu_1-\mu_2|\,|\mu_1+\mu_2| \ge 2\underline\mu\,|\mu_1-\mu_2|,
    \end{equation*}
    it follows that
    \begin{equation*}
        \|\varphi_2(\lambda_1,\mu_1)-\varphi_2(\lambda_2,\mu_2)\|_2^2 \ge (\lambda_1-\lambda_2)^2+4\underline\mu^2(\mu_1-\mu_2)^2 \ge \min\{1,4\underline\mu^2\}\|(\lambda_1,\mu_1)-(\lambda_2,\mu_2)\|_2^2.
    \end{equation*}
    By Theorem~\ref{theorem-ee1}, the objective of reformulation \ref{opt-gg1-r2} is convex under Assumption~\ref{assumption:mu_hat}.\ref{assumption:mu_hat convex}, and $\kappa_2$-strongly convex under Assumption~\ref{assumption:mu_hat}.\ref{assumption:mu_hat strongly convex}. Applying Lemma~\ref{lemma-PLK composition}, we conclude that the original problem \eqref{opt-main problem} satisfies the P\L K condition with exponent $\alpha=1$ and constant
    \begin{equation*}
        \frac{\min\{1,2\underline\mu\}} {2\sqrt{(\bar\lambda-\underline\lambda)^2+(\bar\mu^2-\underline\mu^2)^2}}
    \end{equation*}
    under Assumption~\ref{assumption:mu_hat}.\ref{assumption:mu_hat convex}, and satisfies the P\L K condition with exponent $\alpha=2$ and constant $\kappa_2\min\{1,4\underline\mu^2\}$ under Assumption~\ref{assumption:mu_hat}.\ref{assumption:mu_hat strongly convex}. This completes the proof. \Halmos
\end{proof}

\section{Conclusion}
We study the joint control of arrival and service rates in queueing systems, with the objective of minimizing the long-run expected cost minus revenue. Despite the non-convexity of the objective function in the original decision variables, an appropriate reparameterization yields a convex reformulation. Leveraging this hidden convexity, we show that the original problem satisfies the PŁK condition, which helps explain the observed empirical success of first-order methods under standard assumptions. Our results apply to a broad class of $GI/GI/1$ models, including $\Gamma_k/\Gamma_m/1$ queues with $k, m \ge 1$ and $GI/M/1$ queues with log-concave interarrival times. A key ingredient in our analysis is a new convexity property showing that the expected queue length remains convex after a square-root transformation of the traffic intensity.

Our work opens several directions for future research. A fundamental question is to identify the general characteristics of the interarrival and service time distributions that guarantee hidden convexity. Moving beyond the Gamma family or $GI/M/1$ queues to characterize this property for more general systems would represent a major theoretical advancement. In addition, our findings suggest two paths forward for learning and control. For online learning of static policies, it would be interesting to design algorithms explicitly leveraging the PŁK condition. A more challenging frontier is dynamic pricing, where decisions are state-dependent. It remains an open question whether a similar benign optimization landscape exists in this setting.

% %\THEEndNotes
% \begingroup \parindent 0pt \parskip 0.0ex \def\enotesize{\normalsize} \theendnotes \endgroup

% Appendix here
% Options are (1) APPENDIX (with or without general title) or
%             (2) APPENDICES (if it has more than one unrelated sections)
% Outcomment the appropriate case if necessary
%
% \begin{APPENDIX}{<Title of the Appendix>}
% \end{APPENDIX}
%
%   or
%

% \begin{APPENDICES}
% % \section{<Title of Section B>}
% \end{APPENDICES}

% Acknowledgments here
% \ACKNOWLEDGMENT{We would like to express our sincere gratitude to [acknowledge individuals, organizations, or institutions] for their invaluable contributions to this research. We are also grateful to [mention any additional acknowledgements, such as technical assistance, data providers, or colleagues] for their support and assistance throughout the course of this work.}

\ACKNOWLEDGMENT{We thank Jie Wang for bringing to our attention the empirical observation of global convergence for first-order methods.}

\bibliographystyle{informs2014}
\bibliography{ref}

%%%%%%%%%%%%%%%%%
\end{document}